\documentclass[12pt]{article}

\usepackage{amsmath,amssymb,amsthm,amscd}

\begin{document}

\newcommand{\End}{{\rm{End}\ts}}
\newcommand{\Hom}{{\rm{Hom}}}
\newcommand{\ch}{{\rm{ch}\ts}}
\newcommand{\non}{\nonumber}
\newcommand{\wt}{\widetilde}
\newcommand{\wh}{\widehat}
\newcommand{\ot}{\otimes}
\newcommand{\la}{\lambda}
\newcommand{\La}{\Lambda}
\newcommand{\al}{\alpha}
\newcommand{\be}{\beta}
\newcommand{\ga}{\gamma}
\newcommand{\si}{\sigma}
\newcommand{\vp}{\varphi}
\newcommand{\de}{\delta^{}}
\newcommand{\om}{\omega^{}}
\newcommand{\hra}{\hookrightarrow}
\newcommand{\ve}{\varepsilon}
\newcommand{\ts}{\,}
\newcommand{\qin}{q^{-1}}
\newcommand{\tss}{\hspace{1pt}}
\newcommand{\U}{ {\rm U}}
\newcommand{\Y}{ {\rm Y}}
\newcommand{\CC}{\mathbb{C}\tss}
\newcommand{\SSb}{\mathbb{S}\tss}
\newcommand{\ZZ}{\mathbb{Z}\tss}
\newcommand{\Z}{{\rm Z}}
\newcommand{\Ac}{\mathcal{A}}
\newcommand{\Pc}{\mathcal{P}}
\newcommand{\Qc}{\mathcal{Q}}
\newcommand{\Tc}{\mathcal{T}}
\newcommand{\Sc}{\mathcal{S}}
\newcommand{\Bc}{\mathcal{B}}
\newcommand{\Ec}{\mathcal{E}}
\newcommand{\Hc}{\mathcal{H}}
\newcommand{\Ar}{{\rm A}}
\newcommand{\Ir}{{\rm I}}
\newcommand{\Zr}{{\rm Z}}
\newcommand{\gl}{\mathfrak{gl}}
\newcommand{\Pf}{{\rm Pf}}
\newcommand{\oa}{\mathfrak{o}}
\newcommand{\spa}{\mathfrak{sp}}
\newcommand{\g}{\mathfrak{g}}
\newcommand{\ka}{\mathfrak{k}}
\newcommand{\p}{\mathfrak{p}}
\newcommand{\sll}{\mathfrak{sl}}
\newcommand{\agot}{\mathfrak{a}}
\newcommand{\qdet}{ {\rm qdet}\ts}
\newcommand{\sdet}{ {\rm sdet}\ts}
\newcommand{\Gr}{ {\rm Gr}\tss}
\newcommand{\sgn}{ {\rm sgn}\ts}
\newcommand{\Sym}{\mathfrak S}
\newcommand{\fand}{\quad\text{and}\quad}
\newcommand{\Fand}{\qquad\text{and}\qquad}
\newcommand{\vt}{{\tss|\hspace{-1.5pt}|\tss}}

\renewcommand{\theequation}{\arabic{section}.\arabic{equation}}

\newtheorem{thm}{Theorem}[section]
\newtheorem{lem}[thm]{Lemma}
\newtheorem{prop}[thm]{Proposition}
\newtheorem{cor}[thm]{Corollary}
\newtheorem{conj}[thm]{Conjecture}

\theoremstyle{definition}
\newtheorem{defin}[thm]{Definition}

\theoremstyle{remark}
\newtheorem{remark}[thm]{Remark}
\newtheorem{example}[thm]{Example}

\newcommand{\bth}{\begin{thm}}
\renewcommand{\eth}{\end{thm}}
\newcommand{\bpr}{\begin{prop}}
\newcommand{\epr}{\end{prop}}
\newcommand{\ble}{\begin{lem}}
\newcommand{\ele}{\end{lem}}
\newcommand{\bco}{\begin{cor}}
\newcommand{\eco}{\end{cor}}
\newcommand{\bde}{\begin{defin}}
\newcommand{\ede}{\end{defin}}
\newcommand{\bex}{\begin{example}}
\newcommand{\eex}{\end{example}}
\newcommand{\bre}{\begin{remark}}
\newcommand{\ere}{\end{remark}}
\newcommand{\bcj}{\begin{conj}}
\newcommand{\ecj}{\end{conj}}

\newcommand{\bal}{\begin{aligned}}
\newcommand{\eal}{\end{aligned}}
\newcommand{\beq}{\begin{equation}}
\newcommand{\eeq}{\end{equation}}
\newcommand{\ben}{\begin{equation*}}
\newcommand{\een}{\end{equation*}}

\newcommand{\bpf}{\begin{proof}}
\newcommand{\epf}{\end{proof}}

\def\beql#1{\begin{equation}\label{#1}}

\title{\Large\bf Littlewood--Richardson polynomials}

\author{{\sc A. I. Molev}\\[10mm]
School of Mathematics and Statistics\\
University of Sydney,
NSW 2006, Australia\\
{\tt alexm\hspace{0.09em}@\hspace{0.1em}maths.usyd.edu.au}
}

\date{} 


\maketitle

\vspace{15 mm}

\begin{abstract}
We introduce a family of rings of symmetric functions
depending on an infinite sequence of parameters.
A distinguished basis of such a ring is comprised
by analogues of the Schur functions.
The corresponding structure coefficients are
polynomials in the parameters which we call the
Littlewood--Richardson polynomials.
We give a combinatorial rule for their calculation
by modifying an earlier result of B.~Sagan and the author.
The new rule provides a formula for these
polynomials which is manifestly positive
in the sense of W.~Graham.
We apply this formula for the calculation
of the product of equivariant Schubert classes
on Grassmannians which implies a stability property
of the structure coefficients.
The first manifestly positive
formula for such an expansion
was given by A.~Knutson and T.~Tao
by using combinatorics of puzzles while
the stability property was not apparent
from that formula. We also use
the Littlewood--Richardson polynomials
to describe the multiplication rule in the algebra
of the Casimir elements for the general linear Lie algebra
in the basis of the quantum immanants
constructed by A.~Okounkov and G.~Olshanski.
\end{abstract}

\newpage

\section{Introduction}
\label{sec:int}
\setcounter{equation}{0}

Let $a=(a_i)$, $i\in\ZZ$ be a sequence of variables. Consider
the ring of polynomials $\ZZ[a]$ in the variables $a_i$ with
integer coefficients. Introduce another infinite set of variables
$x=(x_1,x_2,\dots)$ and for each nonnegative integer $n$
denote by $\La_n$ the ring of symmetric
polynomials in $x_1,\dots,x_n$ with coefficients in $\ZZ[a]$.
The ring $\La_n$ is filtered by the usual degrees of polynomials in
$x_1,\dots,x_n$ with the $a_i$ considered to have the zero degree.
The evaluation map
\beql{evhom}
\varphi_n:\La_{n}\to\La_{n-1},\qquad
P(x_1,\dots,x_n)\mapsto
P(x_1,\dots,x_{n-1},a_n)
\eeq
is a homomorphism of filtered rings so that we can define
the inverse limit ring $\La$
by
\beql{ladef}
\La=\lim_{\longleftarrow} \La_n,\qquad n\to\infty,
\eeq
where the limit is taken with respect to the homomorphisms \eqref{evhom}
in the category of filtered rings.
When $a$ is specialized to the sequence of zeros, this reduces to
the usual definition of the ring of symmetric
functions; see e.g. Macdonald~\cite{m:sfh}. In that case,
a distinguished basis of $\La$ is comprised by the Schur functions
$s_{\la}(x)$ parameterized by all partitions $\la$.
The respective analogues of the $s_{\la}(x)$ in the general case are the
{\it double Schur functions\/} $s_{\la}(x\vt a)$ which
form a basis of $\La$ over $\ZZ[a]$. We introduce the
{\it Littlewood--Richardson polynomials\/} $c_{\la\mu}^{\nu}(a)$
as the structure coefficients of the ring $\La$
in the basis of double Schur functions,
\beql{lrca}
s_{\la}(x\vt a)\ts s_{\mu}(x\vt a)=
\sum_{\nu}c_{\la\mu}^{\nu}(a)\ts s_{\nu}(x\vt a).
\eeq
In the specialization $a=(0)$ the polynomials
$c_{\la\mu}^{\nu}(a)$
become the classical Littlewood--Richardson coefficients
$c_{\la\mu}^{\nu}$; see \cite{lr:gc}.
These are remarkable nonnegative integers which occupy a prominent
place in combinatorics, representation theory and geometry;
see e.g. Fulton~\cite{f:yt}, Macdonald~\cite{m:sfh} and
Sagan~\cite{s:sg}.

The main result of this paper is a combinatorial
rule for the calculation
of the Littlewood--Richardson polynomials
which provides a manifestly positive formula
in the sense that $c_{\la\mu}^{\nu}(a)$ is written as
a polynomial in the differences $a_i-a_j$, $i<j$, with
positive integer coefficients.

We consider two applications of the rule.
The results of Knutson and Tao~\cite{kt:pe}
imply that under an appropriate
specialization, the polynomials $c_{\la\mu}^{\nu}(a)$
describe the multiplication rule for
the equivariant Schubert classes
on Grassmannians; see also Fulton~\cite{f:ec} for a more
direct argument. Let $n$ and $N$ be nonnegative
integers with $n\leqslant N$
and let
$\Gr(n,N)$ denote the Grassmannian
of the $n$-dimensional vector subspaces of $\CC^N$.
The torus $T=(\CC^*)^N$ acts naturally on $\Gr(n,N)$.
The equivariant cohomology ring $H^*_T(\Gr(n,N))$
is a module over the polynomial ring
$\ZZ[t_1,\dots,t_N]$ which can be identified with $H^*_T(\{pt\})$,
the equivariant cohomology ring
of a point. This module
has a basis of the equivariant Schubert classes
$\si_{\la}$ parameterized by all diagrams $\la$
contained in the $n\times m$ rectangle, $m=N-n$; see
e.g. \cite{f:yt, f:ec}. Then
\beql{silasimu}
\si_{\la}\ts\si_{\mu}=\sum_{\nu}\ts
d_{\la\mu}^{\ts\nu}\ts\si_{\nu},
\eeq
where $d_{\la\mu}^{\ts\nu}=c_{\la\mu}^{\ts\nu}(a)$
with the sequence $a$ specialized by
\beql{at}
a_{-m+1}=-t_{1},\ts a_{-m+2}=-t_{2},\quad\dots,\quad a_{n}=-t_N,
\eeq
while the remaining parameters $a_i$ are set to zero
(the $t_i$ should be replaced with $y_i$ in the notation
of \cite{kt:pe}). The coefficients $d_{\la\mu}^{\ts\nu}$
are given explicitly as
polynomials in the $t_i-t_j$,
$i>j$, with positive integer coefficients.
This positivity property was established
by Graham~\cite{g:pe}
in the general context of the equivariant Schubert calculus.
The first manifestly positive formula
for the coefficients in the
expansion \eqref{silasimu}
was obtained by Knutson and Tao~\cite{kt:pe}
by using combinatorics of {\it puzzles}.
An earlier rule of Molev and Sagan~\cite{ms:lr}
also calculates $d_{\la\mu}^{\ts\nu}$ but
lacks the explicit positivity property.
Our new rule implies a {\it stability property\/}
of the coefficients $d_{\la\mu}^{\ts\nu}$
(see Corollary~\ref{cor:shub} below).
Even though this property was not pointed out in \cite{kt:pe},
it can be derived directly from the puzzle rule; see also
Fulton~\cite{f:ec} for its geometrical interpretation and
an extension to the equivariant Schubert calculus
on the flag variety.

As another application,
we obtain a rule for the positive integer expansion
of the product of two (virtual) quantum immanants
(or the corresponding higher Capelli operators)
of Okounkov and Olshanski~\cite{o:qi, oo:ss}; cf. \cite{ms:lr}.
The {\it quantum immanants\/} $\SSb_{\la|n}$ are elements
of the center $\Z(\gl_n)$ of the universal enveloping algebra
$\U(\gl_n)$ parameterized by partitions $\la$ with
at most $n$ parts; see \cite{o:qi}.
The elements $\SSb_{\la|n}$ form a basis of $\Z(\gl_n)$
so that we can define the coefficients $f_{\la\mu}^{\ts\nu}$
by the expansion
\ben
\SSb_{\la|n}\ts \SSb_{\mu|n}=\sum_{\nu} f_{\la\mu}^{\ts\nu}\ts
\SSb_{\nu|n}.
\een
Then $f_{\la\mu}^{\ts\nu}=c_{\la\mu}^{\ts\nu}(a)$
for the specialization $a_i=-i$ for $i\in\ZZ$.
As $n\to\infty$ this yields a multiplication rule for the
{\it virtual quantum immanants\/} $\SSb_{\la}$;
see Section~\ref{subsec:hcap} for the definitions.

We define the double Schur function
$s_{\la}(x\vt a)$ as the sequence
of the {\it double Schur polynomials\/}
\beql{dsp}
s_{\la}(x_1,\dots,x_n\vt a),\qquad n=1,2,\dots,
\eeq
which are compatible with respect
to the homomorphisms \eqref{evhom},
\beql{stabil}
\varphi_n:s_{\la}(x_1,\dots,x_n\vt a)\mapsto
s_{\la}(x_1,\dots,x_{n-1}\vt a).
\eeq
The polynomials \eqref{dsp} are closely related to the ``factorial"
or ``double" Schur polynomials $s_{\la}(x|u)$
with $x=(x_1,\dots,x_n)$. The latter were introduced by
Goulden and Greene~\cite{gg:nt} and Macdonald~\cite{m:sf}
as a generalization of the factorial Schur polynomials
of Biedenharn and Louck~\cite{bl:nc, bl:ib},
and they are also a special case of
the double Schubert polynomials of Lascoux and
Sch\"utzenberger; see Lascoux~\cite{l:i}.
We follow Chen, Li and Louck~\cite{cll:fd} and
Fulton~\cite{f:ec} and use
the name ``double Schur polynomials" for the related
polynomials $s_{\la}(x\vt a)$ as well.

In a more detail,
consider a {\it partition\/}
$\la$ which is a sequence
$\la=(\la_1,\dots,\la_n)$ of integers $\la_i$ such that
$\la_1\geqslant\dots\geqslant\la_n\geqslant 0$.
We will identify $\la$ with its {\it diagram\/}
represented graphically as
the array of left justified rows of unit boxes
with $\la_1$ boxes in the top row, $\la_2$ boxes in the second
row, etc. The total number of boxes in $\la$
will be denoted by $|\la|$.
The transposed diagram
$\la'=(\la'_1,\dots,\la'_p)$ is obtained from $\la$
by applying the symmetry
with respect to the main diagonal, so that
$\la'_j$ is the number of boxes in the $j$-th column of $\la$.

Let $u=(u_1,u_2,\dots)$ be a sequence of variables.
The polynomials $s_{\la}(x|u)$ can be defined by
\beql{sxa}
s_{\lambda}(x|u)=\sum_{T}\prod_{\alpha\in\lambda}
(x^{}_{T(\alpha)}-u^{}_{T(\alpha)+c(\alpha)}),
\eeq
where $T$ runs over all semistandard (column-strict)
tableaux of shape $\lambda$
with entries in $\{1,\dots,n\}$,
$T(\alpha)$ is the entry of $T$ in the box
$\alpha\in\lambda$ and $c(\alpha)=j-i$ is the content of
the box
$\alpha=(i,j)$ in row $i$ and
column $j$.

By a {\it reverse $\la$-tableau\/}
$T$ we will mean the tableau obtained by filling in the boxes
of $\la$ with the numbers $1,2,\dots,n$ in such a way that
the entries weakly decrease along the rows and strictly decrease
down the columns. If $\al=(i,j)$ is a box of $\la$
we let $T(\al)=T(i,j)$ denote the entry of $T$ in the box $\al$.
We define the double Schur polynomials $s_{\la}(x\vt a)$ by
\beql{defdouble}
s_{\la}(x\vt a)=\sum_{T}
\prod_{\al\in\la}
(x^{}_{T(\alpha)}-a^{}_{T(\alpha)-c(\al)}),
\eeq
summed over the reverse $\la$-tableaux $T$.
Then we have
\beql{ssassb}
s_{\la}(x\vt a)=s_{\la}(x|u)
\eeq
for the sequences $a$ and $u$ related by
$a_{n-i+1}=u_i$ with $i=1,2,\dots$.
In particular, the polynomial
$s_{\la}(x\vt a)$ only depends on the
variables $a_i$ with
$i\leqslant n$, $i\in\ZZ$.
The relation \eqref{ssassb} is verified easily by
replacing $x_i$ with $x_{n-i+1}$
in \eqref{sxa} for all $i=1,\dots,n$ and using
the fact that $s_{\la}(x|u)$ is a symmetric polynomial in $x$.
The property \eqref{stabil} of the double Schur
polynomials is immediate from their definition.
In the specialization of the sequence $a$ with $a_i=-i$, $i\in\ZZ$,
formula \eqref{defdouble} defines
the shifted Schur polynomials of Okounkov and
Olshanski~\cite{o:qi, oo:ss} in the variables $y_i=x_i+i$.
The use of the {\it reverse\/}
tableaux was significant in their
study of the vanishing and stability properties
of these polynomials and associated central elements
of the universal enveloping algebra for the Lie algebra $\gl_n$;
see also Section~\ref{subsec:hcap} below.

Note that the stability property \eqref{stabil} extends
to the double Schubert polynomials
(and to the equivariant Schubert calculus
on the flag manifold). This follows easily
from the Cauchy formula for the Schubert polynomials
(e.g., put $x_1=y_1$ in \cite[Formula in 2.5.5]{m:sfs}).
In a more general context, this was also pointed out in
\cite{br:sg}.

The double Schur polynomials
$s_{\la}(x\vt a)$ parameterized
by the diagrams $\la$ with at most $n$ rows
form a basis of the ring $\La_n$.
Due to the stability property \eqref{stabil}, the
Littlewood--Richardson polynomials $c_{\la\mu}^{\nu}(a)$
can be defined by the expansion \eqref{lrca},
where $x$ is understood as the set
of variables $x=(x_1,\dots,x_n)$ for any positive integer $n$
such that the diagrams
$\la$, $\mu$ and $\nu$ have at most $n$ rows.
This allows us to work with a finite set of variables
for the determination of the polynomials $c_{\la\mu}^{\nu}(a)$.
For the proof of the main theorem (Theorem~\ref{thm:lrr})
we follow the general approach of \cite{ms:lr},
using the techniques of ``barred" tableaux
and modify the corresponding arguments in order to obtain
manifestly positive polynomials. This is achieved
by imposing a {\it boundness condition\/}
on the barred tableaux.

It was observed by
Goulden and Greene~\cite{gg:nt} and
Macdonald~\cite{m:sf} that $s_{\la}(x|u)$, regarded
as a formal power series
in the infinite sets of variables $x$ and $u$, admits a
``supertableaux" representation. We show that
this representation has its ``finite" counterpart
where $x$ is a finite set of variables.
We derive the corresponding formula
by choosing a certain specialization of the 9th Variation
in \cite{m:sf}.
This representation leads to a
``supertableau" expression for the
Littlewood--Richardson polynomials $c_{\la\mu}^{\nu}(a)$,
although that expression is neither manifestly positive, nor
stable.

After the first version of this paper was completed we have learned
of an independent work of V.~Kreiman~\cite{k:elr},
where a positive equivariant Littlewood--Richardson rule
was given. That rule is equivalent to our Theorem~\ref{thm:lrr}
although the proof in \cite{k:elr} is different. Moreover,
Kreiman's paper also provides a weight-preserving bijection
between the Knutson--Tao puzzles and the
barred tableaux used in Theorem~\ref{thm:lrr}.

\medskip

This work was inspired by Bill Fulton's
lectures~\cite{f:ec}. I am grateful to Bill for
stimulating discussions.

\section{Multiplication rule}
\label{sec:mr}
\setcounter{equation}{0}

Let $R$ denote a sequence of diagrams
\beql{r}
\mu=\rho^{(0)}\to\rho^{(1)}\to
\dots\to\rho^{(l-1)}\to\rho^{(l)}=\nu,
\eeq
where $\rho\to\sigma$ means that $\si$ is obtained from $\rho$
by adding one box.
Let $r_i$ denote the row number of the box added to
the diagram $\rho^{(i-1)}$.
The sequence $r_1r_2\dots r_l$ is called the {\it Yamanouchi symbol\/}
of $R$.
Introduce the ordering on the set of boxes of
a diagram $\la$ by reading them by columns from
left to right and from bottom to top in each column.
We call this the {\it column order\/}.
We shall write $\al\prec \be$ if $\al$ (strictly) precedes $\be$
with respect to the column order. Given a sequence $R$,
construct the set $\Tc(\la,R)$ of {\it barred\/}
reverse $\la$-tableaux
$T$ with entries from $\{1,2,\dots\}$ such that $T$ contains
boxes $\al_1,\dots,\al_l$ with
\ben
\al_1\prec\dots\prec\al_l\Fand T(\al_i)=r_i,\quad
1\leqslant i\leqslant l.
\een
We will distinguish the entries in $\al_1,\dots,\al_l$ by barring
each of them. So, an element of $\Tc(\la,R)$ is a pair consisting
of a reverse $\la$-tableau and a chosen sequence of barred entries
compatible with $R$.
We shall keep the notation $T$ for such a pair.
For example, let $R$ be the sequence
\ben
(3,1)\to(3,2)\to(3,2,1)\to (3,3,1)\to(4,3,1)
\een
so that the Yamanouchi symbol is $2\ts3\ts2\ts1$. Then for
$\la=(5,5,3)$ the following barred $\la$-tableau
belongs to $\Tc(\la,R)$:

\setlength{\unitlength}{0.75em}
\begin{center}
\begin{picture}(18,6.2)

\put(6.7,0.5){$\overline 2$}
\put(6.7,2.5){4}
\put(6.7,4.5){5}
\put(8.7,0.5){1}
\put(8.7,2.5){$\overline 3$}
\put(8.7,4.5){5}
\put(10.7,0.5){1}
\put(10.7,2.5){2}
\put(10.7,4.5){4}
\put(12.7,2.5){1}
\put(12.7,4.5){$\overline 2$}
\put(14.7,2.5){$\overline 1$}
\put(14.7,4.5){2}

\put(6,0){\line(0,1){6}}
\put(8,0){\line(0,1){6}}
\put(10,0){\line(0,1){6}}
\put(12,0){\line(0,1){6}}
\put(14,2){\line(0,1){4}}
\put(16,2){\line(0,1){4}}

\put(6,0){\line(1,0){6}}
\put(6,2){\line(1,0){10}}
\put(6,4){\line(1,0){10}}
\put(6,6){\line(1,0){10}}

\put(16,0){.}

\end{picture}
\end{center}
\setlength{\unitlength}{1pt}

\noindent
For each box $\alpha$ with $\al_i\prec\al\prec\al_{i+1}$,
$0\leqslant i\leqslant l$, set $\rho(\alpha)=\rho^{(i)}$.
The barred entries $\overline r_1,\dots,\overline r_l$
divide the tableau into regions
marked by the elements
of the sequence $R$, as illustrated:

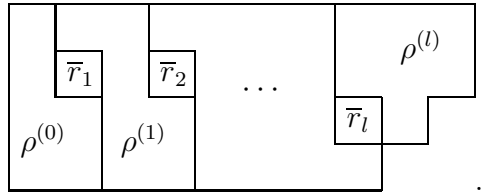
\setlength{\unitlength}{0.75em}
\begin{center}
\begin{picture}(24,8.2)

\put(2,0){\line(0,1){8}}
\put(18,0){\line(0,1){4}}
\put(20,2){\line(0,1){2}}
\put(22,4){\line(0,1){4}}

\put(4,4){\line(0,1){4}}
\put(6,0){\line(0,1){6}}

\put(4,4){\line(1,0){2}}
\put(4,6){\line(1,0){2}}

\put(2.5,1.7){$\rho^{(0)}$}

\put(6.8,1.7){$\rho^{(1)}$}

\put(4.5,4.7){$\overline r_1$}

\put(8,4){\line(0,1){4}}
\put(10,0){\line(0,1){6}}

\put(8,4){\line(1,0){2}}
\put(8,6){\line(1,0){2}}

\put(8.5,4.7){$\overline r_2$}

\put(16,2){\line(0,1){6}}

\put(16,2){\line(1,0){2}}
\put(16,4){\line(1,0){2}}

\put(16.5,2.7){$\overline r_l$}

\put(18.8,5.7){$\rho^{(l)}$}

\put(12,4){$\cdots$}

\put(2,0){\line(1,0){16}}
\put(18,2){\line(1,0){2}}
\put(20,4){\line(1,0){2}}
\put(2,8){\line(1,0){20}}

\put(22,0){.}

\end{picture}
\end{center}
\setlength{\unitlength}{1pt}

\noindent
Finally, a reverse $\la$-tableau $T$ will be called
$\nu$-{\it bounded\/} if
\ben
T(1,j)\leqslant \nu'_j\qquad\text{for all}\quad j=1,\dots,\la_1.
\een
Note that $\nu$-bounded $\la$-tableaux
exist only if $\la\subseteq\nu$.

We are now in a position to state
a rule for the calculation of the
Littlewood-Richardson polynomials
$c_{\la\mu}^{\nu}(a)$ defined by \eqref{lrca}.

\bth\label{thm:lrr}
The polynomial $c_{\la\mu}^{\nu}(a)$ is zero
unless $\mu\subseteq\nu$. If $\mu\subseteq\nu$ then
\beql{lrrule}
c^{\nu}_{\la\mu}(a)=\sum_R\sum_{T}
\prod_{\underset{\scriptstyle T(\alpha)
\text{\ts unbarred}}{\alpha\in\la}}
\Big(\ts a^{}_{\tss T(\alpha)-\rho(\alpha)^{}_{T(\al)}}
-a^{}_{\tss T(\alpha)-c(\alpha)}\ts\Big),
\eeq
summed over all sequences $R$ of the form \eqref{r} and
all $\nu$-bounded
reverse $\la$-tableaux $T\in\Tc(\la,R)$.
Moreover, for each factor occurring
in the formula \eqref{lrrule}
we have $\rho(\al)^{}_{T(\al)}>c(\al)$.
\eth

Before proving the theorem, let us point out some
properties of the Littlewood-Richardson polynomials
which are immediate
from the rule and consider some examples.
The polynomial $c_{\la\mu}^{\nu}(a)$
is zero unless both
diagrams $\la$ and $\mu$ are contained in $\nu$ and
$|\la|+|\mu|\geqslant |\nu|$. In this case
$c_{\la\mu}^{\nu}(a)$ is a homogeneous polynomial
in the $a_i$ of degree $|\la|+|\mu|- |\nu|$.
If $|\la|+|\mu|- |\nu|=0$ then the theorem reproduces
a version of the classical Littlewood-Richardson rule;
see Corollary~\ref{cor:lrclass} below.
Note also that by the definition,
the polynomials have the symmetry
$c_{\la\mu}^{\nu}(a)=c_{\mu\la}^{\nu}(a)$ which is not
apparent from the rule.

\bex\label{ex:calc}
For the product of the double Schur
functions $s_{(2)}(x\vt a)$ and $s_{(2,1)}(x\vt a)$ we have
\ben
\bal
s_{(2)}(x\vt a)\ts s_{(2,1)}(x\vt a)
{}&=s_{(4,1)}(x\vt a)+s_{(3,2)}(x\vt a)+s_{(3,1,1)}(x\vt a)
+s_{(2,2,1)}(x\vt a)\\
{}&+\big(a_{-1}-a_2+a_{-2}-a_0\big)\ts s_{(3,1)}(x\vt a)
+\big(a_{-1}-a_2\big)\ts s_{(2,2)}(x\vt a)\\
{}&+\big(a_{-1}-a_0\big)\ts s_{(2,1,1)}(x\vt a)
+\big(a_{-1}-a_2\big)\ts \big(a_{-1}-a_0\big)
\ts s_{(2,1)}(x\vt a).
\eal
\een
For instance, the coefficient of $s_{(3,1)}(x\vt a)$ is calculated
by the following barred $(2)$-tableaux

\setlength{\unitlength}{0.75em}
\begin{center}
\begin{picture}(27,2.2)

\put(2,0){\line(0,1){2}}
\put(4,0){\line(0,1){2}}
\put(6,0){\line(0,1){2}}

\put(2,0){\line(1,0){4}}
\put(2,2){\line(1,0){4}}

\put(12,0){\line(0,1){2}}
\put(14,0){\line(0,1){2}}
\put(16,0){\line(0,1){2}}

\put(12,0){\line(1,0){4}}
\put(12,2){\line(1,0){4}}

\put(22,0){\line(0,1){2}}
\put(24,0){\line(0,1){2}}
\put(26,0){\line(0,1){2}}

\put(22,0){\line(1,0){4}}
\put(22,2){\line(1,0){4}}

\put(2.7,0.5){1}
\put(4.7,0.5){$\overline 1$}

\put(12.7,0.5){$\overline 1$}
\put(14.7,0.5){1}

\put(22.7,0.5){2}
\put(24.7,0.5){$\overline 1$}

\end{picture}
\end{center}
\setlength{\unitlength}{1pt}

\noindent
compatible with the sequence $(2,1)\to(3,1)$. They
contribute respectively $a_{-1}-a_1$, $a_{-2}-a_0$,
$a_1-a_{2}$ which sums up to the coefficient
$a_{-1}-a_2+a_{-2}-a_0$. Alternatively,
using the symmetry $c_{\la\mu}^{\nu}(a)=c_{\mu\la}^{\nu}(a)$
we can calculate the coefficient of
$s_{(3,1)}(x\vt a)$ by considering
the barred $(2,1)$-tableaux

\setlength{\unitlength}{0.75em}
\begin{center}
\begin{picture}(27,4.2)

\put(2,0){\line(0,1){4}}
\put(4,0){\line(0,1){4}}
\put(6,2){\line(0,1){2}}

\put(2,0){\line(1,0){2}}
\put(2,2){\line(1,0){4}}
\put(2,4){\line(1,0){4}}

\put(22,0){\line(0,1){4}}
\put(24,0){\line(0,1){4}}
\put(26,2){\line(0,1){2}}

\put(22,0){\line(1,0){2}}
\put(22,2){\line(1,0){4}}
\put(22,4){\line(1,0){4}}

\put(2.7,0.5){$\overline 1$}
\put(2.7,2.5){$\overline 2$}
\put(4.7,2.5){1}

\put(22.7,0.5){1}
\put(22.7,2.5){$\overline 2$}
\put(24.7,2.5){$\overline 1$}

\end{picture}
\end{center}
\setlength{\unitlength}{1pt}

\noindent
compatible with the sequences $(2)\to (3)\to (3,1)$
and $(2)\to (2,1)\to (3,1)$, respectively.
Their contributions to the coefficient are
$a_{-2}-a_0$ and $a_{-1}-a_2$.
\eex

\medskip
\bex\label{ex:calctwo}
For the calculation of $c_{(4,2,1)\tss(2,2)}^{(5,2,2)}(a)$
take $\la=(4,2,1)$, $\mu=(2,2)$ and $\nu=(5,2,2)$.
We have ten sequences $R$ of the form \eqref{r}
but the set $\Tc(\la,R)$ contains $\nu$-bounded
tableaux only for three of them.
For the sequence $R_1$ with the Yamanouchi symbol
$1\ts3\ts3\ts1\ts1$, the set $\Tc(\la,R_1)$
contains two bounded barred tableaux

\setlength{\unitlength}{0.75em}
\begin{center}
\begin{picture}(31,6.2)

\put(2,0){\line(0,1){6}}
\put(4,0){\line(0,1){6}}
\put(6,2){\line(0,1){4}}
\put(8,4){\line(0,1){2}}
\put(10,4){\line(0,1){2}}

\put(2,0){\line(1,0){2}}
\put(2,2){\line(1,0){4}}
\put(2,4){\line(1,0){8}}
\put(2,6){\line(1,0){8}}

\put(22,0){\line(0,1){6}}
\put(24,0){\line(0,1){6}}
\put(26,2){\line(0,1){4}}
\put(28,4){\line(0,1){2}}
\put(30,4){\line(0,1){2}}

\put(22,0){\line(1,0){2}}
\put(22,2){\line(1,0){4}}
\put(22,4){\line(1,0){8}}
\put(22,6){\line(1,0){8}}

\put(2.7,0.5){$\overline 1$}
\put(2.7,2.5){2}
\put(2.7,4.5){$\overline 3$}
\put(4.7,2.5){2}
\put(4.7,4.5){$\overline 3$}
\put(6.7,4.5){$\overline 1$}
\put(8.7,4.5){$\overline 1$}

\put(22.7,0.5){$\overline 1$}
\put(22.7,2.5){2}
\put(22.7,4.5){$\overline 3$}
\put(24.7,2.5){1}
\put(24.7,4.5){$\overline 3$}
\put(26.7,4.5){$\overline 1$}
\put(28.7,4.5){$\overline 1$}

\end{picture}
\end{center}
\setlength{\unitlength}{1pt}

\noindent
whose contributions to the
Littlewood--Richardson polynomial are $(a_0-a_3)(a_0-a_2)$ and
$(a_0-a_3)(a_{-2}-a_1)$, respectively.
For the sequence $R_2$ with the Yamanouchi symbol
$1\ts3\ts1\ts3\ts1$, the set $\Tc(\la,R_2)$
contains the bounded tableaux

\setlength{\unitlength}{0.75em}
\begin{center}
\begin{picture}(31,6.2)

\put(2,0){\line(0,1){6}}
\put(4,0){\line(0,1){6}}
\put(6,2){\line(0,1){4}}
\put(8,4){\line(0,1){2}}
\put(10,4){\line(0,1){2}}

\put(2,0){\line(1,0){2}}
\put(2,2){\line(1,0){4}}
\put(2,4){\line(1,0){8}}
\put(2,6){\line(1,0){8}}

\put(22,0){\line(0,1){6}}
\put(24,0){\line(0,1){6}}
\put(26,2){\line(0,1){4}}
\put(28,4){\line(0,1){2}}
\put(30,4){\line(0,1){2}}

\put(22,0){\line(1,0){2}}
\put(22,2){\line(1,0){4}}
\put(22,4){\line(1,0){8}}
\put(22,6){\line(1,0){8}}

\put(2.7,0.5){$\overline 1$}
\put(2.7,2.5){2}
\put(2.7,4.5){$\overline 3$}
\put(4.7,2.5){$\overline 1$}
\put(4.7,4.5){$\overline 3$}
\put(6.7,4.5){$\overline 1$}
\put(8.7,4.5){1}

\put(22.7,0.5){$\overline 1$}
\put(22.7,2.5){2}
\put(22.7,4.5){$\overline 3$}
\put(24.7,2.5){$\overline 1$}
\put(24.7,4.5){$\overline 3$}
\put(26.7,4.5){1}
\put(28.7,4.5){$\overline 1$}

\end{picture}
\end{center}
\setlength{\unitlength}{1pt}

\noindent
with the respective contributions $(a_0-a_3)(a_{-4}-a_{-2})$ and
$(a_0-a_3)(a_{-3}-a_{-1})$.
For the sequence $R_3$ with the Yamanouchi symbol
$3\ts1\ts3\ts1\ts1$, the set $\Tc(\la,R_3)$
contains the only bounded tableau

\setlength{\unitlength}{0.75em}
\begin{center}
\begin{picture}(30,6.2)

\put(12,0){\line(0,1){6}}
\put(14,0){\line(0,1){6}}
\put(16,2){\line(0,1){4}}
\put(18,4){\line(0,1){2}}
\put(20,4){\line(0,1){2}}

\put(12,0){\line(1,0){2}}
\put(12,2){\line(1,0){4}}
\put(12,4){\line(1,0){8}}
\put(12,6){\line(1,0){8}}

\put(12.7,0.5){1}
\put(12.7,2.5){2}
\put(12.7,4.5){$\overline 3$}
\put(14.7,2.5){$\overline 1$}
\put(14.7,4.5){$\overline 3$}
\put(16.7,4.5){$\overline 1$}
\put(18.7,4.5){$\overline 1$}

\end{picture}
\end{center}
\setlength{\unitlength}{1pt}

\noindent
with the contribution $(a_{-1}-a_3)(a_0-a_3)$.
Hence,
\ben
c_{(4,2,1)\tss(2,2)}^{(5,2,2)}(a)=(a_0-a_3)\ts
(a_{-4}+a_{-3}+a_0-a_1-a_2-a_3).
\een
Taking $\la=(2,2)$, $\mu=(4,2,1)$ and $\nu=(5,2,2)$
we get two sequences with the
Yamanouchi symbols
$1\ts3$ and $3\ts1$.
The corresponding sets $\Tc(\la,R)$ consist of
five and four bounded barred tableaux, respectively,
thus leading to a slightly longer calculation.
\qed
\eex

\bpf[Proof of Theorem~\ref{thm:lrr}]
We present the proof as a sequence of lemmas.
Due to the stability property \eqref{stabil}, we may
(and will) work
with a finite set of variables $x=(x_1,\dots,x_n)$.
Accordingly, possible entries of the tableaux are now
elements of the set $\{1,\dots,n\}$.
Introduce another sequence of variables $b=(b_i)$, $i\in\ZZ$,
and define the Littlewood--Richardson type
coefficients $c_{\la\mu}^{\nu}(a,b)$ by the expansion
\beql{lrcab}
s_{\la}(x\vt b)\ts s_{\mu}(x\vt a)=
\sum_{\nu}c_{\la\mu}^{\nu}(a,b)\ts s_{\nu}(x\vt a).
\eeq

\ble\label{lem:msab}
The coefficient $c_{\la\mu}^{\nu}(a,b)$ is zero
unless $\mu\subseteq\nu$. If $\mu\subseteq\nu$ then
\beql{lrruleab}
c^{\nu}_{\la\mu}(a,b)=\sum_R\sum_{T}
\prod_{\underset{\scriptstyle T(\alpha)
\text{\ts unbarred}}{\alpha\in\la}}
\Big(\ts a^{}_{\tss T(\alpha)-\rho(\alpha)^{}_{T(\al)}}
-b^{}_{\tss T(\alpha)-c(\alpha)}\ts\Big),
\eeq
summed over all sequences $R$ of the form \eqref{r} and
all reverse $\la$-tableaux $T\in\Tc(\la,R)$.
\ele

\bpf
This is essentially a reformulation of the main result
of \cite{ms:lr} (Theorem~3.1). Note that
the summation in \eqref{lrruleab} is taken over {\it all\/}
barred tableaux $T\in\Tc(\la,R)$ (not just over the
$\nu$-bounded ones as in \eqref{lrrule}).
Rather than repeating the whole argument of \cite{ms:lr},
we only sketch the main steps of the proof and indicate
the necessary changes to be made. We refer the reader
to \cite{ms:lr} for the details.

We assume that all diagrams here have at most $n$ rows.
If $\rho=(\rho_1,\dots,\rho_n)$ is a such diagram,
we set
\ben
a_{\rho}=(a_{1-\rho_1},\dots,a_{n-\rho_n})
\Fand
|a_{\rho}|=a_{1-\rho_1}+\dots+a_{n-\rho_n}.
\een
Under the correspondence \eqref{ssassb} we have
$a_{\rho}=u_{\rho}=(u_{\rho_1+n},\dots,u_{\rho_n+1})$,
the latter notation was used in \cite{ms:lr}.

The starting point is the Vanishing Theorem
of \cite{o:qi} whose proof was also
reproduced in \cite{ms:lr}. By that theorem,
\ben
s_{\la}(a_{\rho}\vt a)=0\qquad\text{unless}
\quad \la\subseteq\rho,
\een
and
\ben
s_{\la}(a_{\la}\vt a)=\prod_{(i,j)\in\la}
\big(a^{}_{i-\la_i}-a^{}_{\la'_j-j+1}\big).
\een
The first claim of the lemma follows
from the Vanishing Theorem which also implies
\ben
c_{\la\mu}^{\mu}(a,b)=s_{\la}(a_{\mu}\vt b).
\een
This proves \eqref{lrruleab} for the case $\nu=\mu$.
Now we suppose that $|\nu|-|\mu|\geqslant 1$
and proceed by induction on $|\nu|-|\mu|$.
The induction step is based on the recurrence relation
\beql{recrel}
c_{\la\mu}^{\nu}(a,b)=\frac{1}{|a_{\nu}|-|a_{\mu}|}
\Bigg(\sum_{\mu\to\mu^+}c_{\la\mu^+}^{\nu}(a,b)
-\sum_{\nu^-\to\nu}c_{\la\mu}^{\nu^-}(a,b)\Bigg)
\eeq
which was proved in \cite[Proposition~3.4]{ms:lr};
see also \cite{kt:pe}.
Suppose that the diagram $\nu$
is obtained from $\mu$ by adding one box in row $r$.
Then
\beql{stepone}
c_{\la\mu}^{\nu}(a,b)=
\frac{s_{\la}(a_{\nu}\vt b)-s_{\la}(a_{\mu}\vt b)}
{(a_{\nu})_r-(a_{\mu})_r}.
\eeq
Now use the definition \eqref{defdouble} of the double Schur
polynomials. Since the $n$-tuples $a_{\nu}$ and $a_{\mu}$
only differ at the $r$-th component,
the ratio on the right hand side of
\eqref{stepone} can be expanded by taking into
account the entries $r$ of the reverse $\la$-tableaux $T$.
We need the
following formula, where we are thinking of
$y=(a_{\nu})_{r}$, $z=(a_{\mu})_r$ and
$m_i=b_{T(\al)-c(\al)}$ as $\al$ runs over
the boxes of $T$ with $T(\al)=r$ in column order:
\ben
\frac{\prod_{i=1}^k (y-m_i)-\prod_{i=1}^k (z-m_i)}{y-z}=
\sum_{j=1}^k (z-m_1)\dots(z-m_{j-1})(y-m_{j+1})\dots(y-m_k).
\een
The right hand side of \eqref{stepone}
can now be interpreted as the right hand side of \eqref{lrruleab},
where $R$ is the only sequence $\mu\to\nu$ and the
sum is taken over the reverse $\la$-tableaux $T$
with one barred entry $r$, as illustrated:

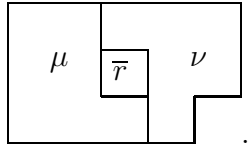
\setlength{\unitlength}{0.75em}
\begin{center}
\begin{picture}(13,6.2)

\put(2,0){\line(0,1){6}}
\put(10,0){\line(0,1){2}}
\put(12,2){\line(0,1){4}}

\put(3.8,3.3){$\mu$}

\put(6,2){\line(0,1){4}}
\put(8,0){\line(0,1){4}}

\put(6,2){\line(1,0){2}}
\put(6,4){\line(1,0){2}}

\put(6.5,2.7){$\overline r$}

\put(9.8,3.3){$\nu$}

\put(2,0){\line(1,0){8}}
\put(10,2){\line(1,0){2}}
\put(2,6){\line(1,0){10}}

\put(12,0){.}

\end{picture}
\end{center}
\setlength{\unitlength}{1pt}

\noindent
Here $\rho(\al)=\mu$ for all boxes $\al$ preceding
the box occupied by the barred $r$, and $\rho(\al)=\nu$
for all boxes $\al$ which follow that box in column order.
Note that the variables $y$ and $z$ are now
swapped on the right
hand side of the above expansion, as compared to \cite{ms:lr}
(this does not change the polynomial due to the symmetry
in $y$ and $z$). Consequently, the column order used
in \cite{ms:lr} is the opposite
to the order on the boxes of $\la$
we use here.

We can represent the above calculation of $c_{\la\mu}^{\nu}(a)$
by the ``diagrammatic" relation

\setlength{\unitlength}{0.75em}
\begin{center}
\begin{picture}(44.5,6.2)

\put(0,3){$\Big(|a_{\nu}|-|a_{\mu}|\Big)$}

\put(8,0){\line(0,1){6}}
\put(16,0){\line(0,1){2}}
\put(18,2){\line(0,1){4}}

\put(9.8,3.3){$\mu$}

\put(12,2){\line(0,1){4}}
\put(14,0){\line(0,1){4}}

\put(12,2){\line(1,0){2}}
\put(12,4){\line(1,0){2}}

\put(12.5,2.7){$\overline r$}

\put(15.8,3.3){$\nu$}

\put(8,0){\line(1,0){8}}
\put(16,2){\line(1,0){2}}
\put(8,6){\line(1,0){10}}

\put(19,3){$=$}

\put(21,0){\line(0,1){6}}
\put(29,0){\line(0,1){2}}
\put(31,2){\line(0,1){4}}

\put(24.8,3.3){$\nu$}

\put(21,0){\line(1,0){8}}
\put(29,2){\line(1,0){2}}
\put(21,6){\line(1,0){10}}

\put(32,3){$-$}

\put(34,0){\line(0,1){6}}
\put(42,0){\line(0,1){2}}
\put(44,2){\line(0,1){4}}

\put(37.8,3.3){$\mu$}

\put(34,0){\line(1,0){8}}
\put(42,2){\line(1,0){2}}
\put(34,6){\line(1,0){10}}

\put(44,0){.}

\end{picture}
\end{center}
\setlength{\unitlength}{1pt}

\noindent
Consider now the next case where $|\nu|-|\mu|=2$
and apply the recurrence relation \eqref{recrel}.
We have three subcases: the diagram $\nu$ is
obtained from $\mu$ by adding two boxes in different rows and
columns; by adding two boxes in the same row;
or by adding two boxes in the same column.
The first two subcases are dealt with in a way similar
to the case $|\nu|-|\mu|=1$. An additional care is needed
for the third subcase where we suppose that $\nu$
is obtained from $\mu$ by adding the boxes in rows $r$ and $r+1$.
Denote by $\rho$ the diagram obtained from $\mu$ by adding
the box in row $r$. Then \eqref{recrel} gives
\ben
c_{\la\mu}^{\nu}(a,b)=
\frac{c_{\la\rho}^{\nu}(a,b)-c_{\la\mu}^{\rho}(a,b)}
{|a_{\nu}|-|a_{\mu}|}.
\een
Set $s=r+1$. Exactly as in the case $|\nu|-|\mu|=1$,
we have the following diagrammatic relations:

\setlength{\unitlength}{0.75em}
\begin{center}
\begin{picture}(44.5,6.2)

\put(0,3){$\Big(|a_{\rho}|-|a_{\mu}|\Big)$}

\put(8,0){\line(0,1){6}}
\put(16,0){\line(0,1){5}}
\put(18,2){\line(0,1){4}}

\put(8.5,3.3){$\mu$}

\put(10,1){\line(0,1){5}}
\put(12,0){\line(0,1){3}}
\put(10,1){\line(1,0){2}}
\put(10,3){\line(1,0){2}}

\put(10.5,1.7){$\overline r$}

\put(14,3){\line(0,1){3}}
\put(14,3){\line(1,0){2}}
\put(14,5){\line(1,0){2}}

\put(14.5,3.7){$\overline {s}$}

\put(12.5,3.3){$\rho$}

\put(16.8,3.3){$\nu$}

\put(8,0){\line(1,0){8}}
\put(16,2){\line(1,0){2}}
\put(8,6){\line(1,0){10}}

\put(19,3){$=$}

\put(22.8,3.3){$\rho$}

\put(25,2){\line(0,1){4}}
\put(27,0){\line(0,1){4}}
\put(25,2){\line(1,0){2}}
\put(25,4){\line(1,0){2}}

\put(25.5,2.7){$\overline s$}

\put(21,0){\line(0,1){6}}
\put(29,0){\line(0,1){2}}
\put(31,2){\line(0,1){4}}

\put(28.8,3.3){$\nu$}

\put(21,0){\line(1,0){8}}
\put(29,2){\line(1,0){2}}
\put(21,6){\line(1,0){10}}

\put(32,3){$-$}

\put(35.8,3.3){$\mu$}

\put(38,2){\line(0,1){4}}
\put(40,0){\line(0,1){4}}
\put(38,2){\line(1,0){2}}
\put(38,4){\line(1,0){2}}

\put(38.5,2.7){$\overline s$}

\put(34,0){\line(0,1){6}}
\put(42,0){\line(0,1){2}}
\put(44,2){\line(0,1){4}}

\put(41.8,3.3){$\nu$}

\put(34,0){\line(1,0){8}}
\put(42,2){\line(1,0){2}}
\put(34,6){\line(1,0){10}}

\end{picture}
\end{center}
\setlength{\unitlength}{1pt}

\noindent
and

\setlength{\unitlength}{0.75em}
\begin{center}
\begin{picture}(44.5,6.2)

\put(0,3){$\Big(|a_{\nu}|-|a_{\rho}|\Big)$}

\put(8,0){\line(0,1){6}}
\put(16,0){\line(0,1){5}}
\put(18,2){\line(0,1){4}}

\put(8.5,3.3){$\mu$}

\put(10,1){\line(0,1){5}}
\put(12,0){\line(0,1){3}}
\put(10,1){\line(1,0){2}}
\put(10,3){\line(1,0){2}}

\put(10.5,1.7){$\overline r$}

\put(14,3){\line(0,1){3}}
\put(14,3){\line(1,0){2}}
\put(14,5){\line(1,0){2}}

\put(14.5,3.7){$\overline {s}$}

\put(12.5,3.3){$\rho$}

\put(16.8,3.3){$\nu$}

\put(8,0){\line(1,0){8}}
\put(16,2){\line(1,0){2}}
\put(8,6){\line(1,0){10}}

\put(19,3){$=$}

\put(22.8,3.3){$\mu$}

\put(25,2){\line(0,1){4}}
\put(27,0){\line(0,1){4}}
\put(25,2){\line(1,0){2}}
\put(25,4){\line(1,0){2}}

\put(25.5,2.7){$\overline r$}

\put(21,0){\line(0,1){6}}
\put(29,0){\line(0,1){2}}
\put(31,2){\line(0,1){4}}

\put(28.8,3.3){$\nu$}

\put(21,0){\line(1,0){8}}
\put(29,2){\line(1,0){2}}
\put(21,6){\line(1,0){10}}

\put(32,3){$-$}

\put(35.8,3.3){$\mu$}

\put(38,2){\line(0,1){4}}
\put(40,0){\line(0,1){4}}
\put(38,2){\line(1,0){2}}
\put(38,4){\line(1,0){2}}

\put(38.5,2.7){$\overline r$}

\put(34,0){\line(0,1){6}}
\put(42,0){\line(0,1){2}}
\put(44,2){\line(0,1){4}}

\put(41.8,3.3){$\rho$}

\put(34,0){\line(1,0){8}}
\put(42,2){\line(1,0){2}}
\put(34,6){\line(1,0){10}}

\put(44,0){.}

\end{picture}
\end{center}
\setlength{\unitlength}{1pt}

\noindent
Hence, the desired formula for $c_{\la\mu}^{\nu}(a,b)$
will follow if we prove the
relation

\setlength{\unitlength}{0.75em}
\begin{center}
\begin{picture}(24.5,6.2)

\put(0,0){\line(0,1){6}}
\put(8,0){\line(0,1){2}}
\put(10,2){\line(0,1){4}}

\put(1.8,3.3){$\mu$}

\put(4,2){\line(0,1){4}}
\put(6,0){\line(0,1){4}}
\put(4,2){\line(1,0){2}}
\put(4,4){\line(1,0){2}}

\put(4.5,2.7){$\overline r$}

\put(7.8,3.3){$\nu$}

\put(0,0){\line(1,0){8}}
\put(8,2){\line(1,0){2}}
\put(0,6){\line(1,0){10}}

\put(11,3){$=$}

\put(13,0){\line(0,1){6}}
\put(21,0){\line(0,1){2}}
\put(23,2){\line(0,1){4}}

\put(14.8,3.3){$\mu$}

\put(17,2){\line(0,1){4}}
\put(19,0){\line(0,1){4}}
\put(17,2){\line(1,0){2}}
\put(17,4){\line(1,0){2}}

\put(17.5,2.7){$\overline s$}

\put(20.8,3.3){$\nu$}

\put(13,0){\line(1,0){8}}
\put(21,2){\line(1,0){2}}
\put(13,6){\line(1,0){10}}

\put(23,0){.}

\end{picture}
\end{center}
\setlength{\unitlength}{1pt}

\noindent
We construct a weight-preserving
bijection between the barred reverse $\la$-tableaux
which are represented by the left and right hand sides
of this diagrammatic relation.
Here the weight is the product on the right hand side
of \eqref{lrruleab} corresponding to a barred tableau.
Let such a tableau
with a barred entry $r$ in the box $(i,j)$ be given.
Suppose first that the box $(i-1,j)$ belongs
to the diagram and it is occupied by $s=r+1$.
Then the image of the tableau under the map
is the same tableau but the entry $T(i,j)=r$ is now unbarred
while $T(i-1,j)=r+1$ is barred. Since
\ben
(a_{\nu})_{r+1}=(a_{\mu})_r\Fand
T(i-1,j)-c(i-1,j)=T(i,j)-c(i,j),
\een
the weights of the
tableaux are preserved under the map.

Suppose now
that the entry in the box $(i-1,j)$ is greater than $r+1$,
or this box is outside the diagram. Consider
all entries $r$ in the row $i$
to the left of the box $(i,j)$ and suppose that they
occupy the boxes $(i,j-m),\ts(i,j-m+1),\dots,(i,j-1)$.
Then the image of the tableau under the map
is the tableau obtained by
replacing the entries in each of the
boxes $(i,j-m),\dots,(i,j)$ with $s=r+1$ and barring
the entry in the box $(i,j-m)$. The weights of the
tableaux are again preserved.

The inverse map is described in a similar way.
This gives the desired weight-preserving bijection.
The general argument uses similar
calculations with the barred diagrams
and a similar bijection described in \cite{ms:lr}.
\epf

\bre\label{rem:ms}
(i)\quad A cohomological interpretation of
the coefficients $c_{\la\mu}^{\nu}(a,b)$
and their puzzle computation can be found in \cite{kt:pe}.
\par
(ii)\quad The definition \eqref{lrcab} of the coefficients
$c_{\la\mu}^{\nu}(a,b)$ can be extended to
the case where $\la$ is a skew diagram.
Lemma~\ref{lem:msab} and its proof remain valid; see \cite{ms:lr}.
\par
(iii)\quad In contrast with
the Littlewood--Richardson polynomials $c_{\la\mu}^{\nu}(a)$,
the coefficients $c_{\la\mu}^{\nu}(a,b)$ do not have the
stability property as they depend on $n$.
\qed
\ere

\medskip

Lemma~\ref{lem:msab} implies that the
Littlewood--Richardson polynomials
can be calculated by \eqref{lrruleab} with $b=a$, that is,
$c_{\la\mu}^{\nu}(a)=c_{\la\mu}^{\nu}(a,a)$.
Our strategy now is to show that (unlike the formula
of Theorem~3.1 in \cite{ms:lr}), the formula \eqref{lrruleab}
(with $b=a$)
is ``nonnegative" in the sense that
all nonzero products which occur in the formula
are polynomials in the $a_i-a_j$ with $i<j$. Then we
demonstrate that the $\nu$-boundness condition
serves to eliminate the unwanted zero terms.

\ble\label{lem:nonneg}
Let $R$ be a sequence of the form \eqref{r} and
let $T\in\Tc(\la,R)$. Suppose that
\beql{prnz}
\prod_{\underset{\scriptstyle T(\alpha)\text{\ts unbarred}}{\alpha\in\la}}
\Big(\ts a^{}_{\tss T(\alpha)-\rho(\alpha)^{}_{T(\al)}}
-a^{}_{\tss T(\alpha)-c(\alpha)}\ts\Big)\ne 0.
\eeq
Then $\rho(\alpha)^{}_{T(\al)}> c(\al)$ for
all $\al\in\la$ with unbarred $T(\alpha)$.
\ele

\bpf
Suppose on the contrary that there exists a box $\al=(i,j)$ with
an unbarred $T(i,j)$ and the condition
$\rho(i,j)_{T(i,j)}<j-i$; the equality $\rho(i,j)_{T(i,j)}=j-i$
is excluded
since this would violate \eqref{prnz}.
Choose such a box with the minimum
possible value of $j$. If all the entries $T(i,1),\dots,T(i,j-1)$
of $T$ are barred then $\rho(i,j)$ is obtained from $\mu$
by adding boxes in rows $T(i,1)\geqslant\dots\geqslant T(i,j-1)$
and, possibly, by adding other boxes. Since
$T(i,j-1)\geqslant T(i,j)$, we have $\rho(i,j)_{T(i,j)}\geqslant j-1$,
a contradiction. So, at least one of the entries
$T(i,1),\dots,T(i,j-1)$ must be unbarred. Take
such an unbarred entry $T(i,k)$ which is the closest to $T(i,j)$,
that is, all entries $T(i,k+1),\dots,T(i,j-1)$ are barred.
Then $\rho(i,j)$ is obtained from $\rho(i,k)$
by adding boxes in rows $T(i,k+1)\geqslant\dots\geqslant T(i,j-1)$
and, possibly, by adding other boxes. Hence,
\ben
\rho(i,j)_{T(i,j)}\geqslant \rho(i,k)_{T(i,k)}+j-k-1
\een
which implies $\rho(i,k)_{T(i,k)}<k-i+1$. However, if
$\rho(i,k)_{T(i,k)}=k-i$ then the factor in \eqref{prnz}
corresponding to $\al=(i,k)$ is zero, which is impossible.
Therefore $\rho(i,k)_{T(i,k)}<k-i$ which contradicts the choice
of $j$.
\epf

\ble\label{lem:necc}
Suppose that $R$ is a sequence of the form \eqref{r} and
$T\in\Tc(\la,R)$. If \eqref{prnz} holds then
$T$ is $\nu$-bounded.
\ele

\bpf
By Lemma~\ref{lem:nonneg}, for all unbarred entries
$T(1,k)$ of the first row of the tableau $T$
we have $\rho(1,k)_{T(1,k)}\geqslant k$.
This implies $\nu_{T(1,k)}\geqslant k$.
If the entry
$T(1,j)$ is barred then $\rho(1,k)_{T(1,k)}\geqslant k$
for the nearest unbarred entry $T(1,k)$ on its left
(if it exists).
Then $\nu$ is obtained from $\rho(1,k)$
by adding boxes in rows $T(1,k+1)\geqslant\dots\geqslant T(1,j)$
and, possibly, by adding other boxes.
This implies $\nu_{T(1,j)}\geqslant j$. Thus, this
inequality holds for all $j=1,\dots,\la_1$.
This is equivalent to the $\nu$-boundness of $T$.
\epf

\ble\label{lem:admsuff}
Suppose that $R$ is a sequence of the form \eqref{r} and
$T\in\Tc(\la,R)$ is $\nu$-bounded.
Then $\rho(\alpha)^{}_{T(\al)}> c(\al)$ for
all $\al\in\la$ with unbarred $T(\alpha)$.
\ele

\bpf
We argue by contradiction. Taking into account
Lemma~\ref{lem:nonneg}, we find that for some
$\al=(i,j)$ with unbarred $T(\alpha)$ we have
$\rho(i,j)_{T(i,j)}=j-i$. Set $t=T(i,j)$ and consider
all barred entries of $T$ (assuming for now they exist)
which are equal to $t$
and occur to the right of the column $j$.
Since $T$ is a reverse tableau, these entries $\bar t$
can only occur in rows $1,2,\dots,i$.
Let $(r,k)$ be the box with the maximum column number $k$
containing $\bar t$.
Then the total number of such entries $\bar t$
does not exceed $k-j$.
This implies that the number of boxes $\nu_t$
in row $t$ of $\nu$
does not exceed $\rho(i,j)_t+k-j=k-i$. Hence,
$\nu'_k\leqslant t-1$. On the other hand,
by the $\nu$-boundness of $T$ we have
$t=T(r,k)\leqslant T(1,k)\leqslant \nu'_k$, a contradiction.

If none of the boxes to the right of the column $j$
contains $\bar t$ then $\nu_t=\rho(i,j)_t=j-i$.
However, by the assumption, $\nu_t\geqslant\nu_{T(1,j)}\geqslant j$,
a contradiction.
\epf

This completes the proof of the theorem.
\epf

By the {\it column word\/} of a tableau $T$ we will mean
the sequence of all entries of $T$ written
in the column order.

\bco\label{cor:lrclass}
Suppose that $|\nu|=|\la|+|\mu|$.
The Littlewood--Richardson coefficient
$c_{\la\mu}^{\nu}$ equals the number of $\nu$-bounded
reverse $\la$-tableaux $T$ whose column word coincides
with the Yamanouchi symbol of a certain sequence $R$
of the form \eqref{r}.
\qed
\eco

This can be shown to be equivalent to
a well-known version of the Littlewood--Richardson rule.
Corollary~\ref{cor:lrclass} also holds
with the $\nu$-boundness condition dropped;
see Lemma~\ref{lem:necc}.
By the corollary, $c_{\la\mu}^{\nu}$
counts the cardinality of the intersection of two
finite sets: the set of column words of
$\nu$-bounded reverse $\la$-tableaux and the set of
Yamanouchi symbols of the sequences of the form \eqref{r}.

\bre\label{rem:lrr}
Due to \eqref{ssassb}, the multiplication rule for
the polynomials $s_{\la}(x|u)$ is obtained
from Theorem~\ref{thm:lrr} by replacing $a_i$ with $u_{n-i+1}$
for each $i$.
The corresponding coefficients are polynomials
in the $u_i-u_j$, $i>j$, with positive integer
coefficients.
\qed
\ere

\bco\label{cor:finn}
Suppose that the
polynomials $c_{\la\mu}^{\nu}(a)$ are defined by the expansion
\eqref{lrca} with $x=(x_1,\dots,x_n)$. Then
$c_{\la\mu}^{\nu}(a)$ is independent
of $n$ as soon as $n \geqslant \nu^{\tss\prime}_1$.
Moreover, if $n < \nu^{\tss\prime}_1$ then $c_{\la\mu}^{\nu}(a)=0$.
\eco

\bpf
This follows from the boundness condition on the reverse
tableaux.
\epf

\section{Applications}
\label{sec:appl}
\setcounter{equation}{0}

\subsection{Equivariant Schubert calculus
on the Grassmannian}

As in the Introduction, consider the
equivariant cohomology ring $H^*_T(\Gr(n,N))$
as a module over $\ZZ[t_1,\dots,t_N]$.
Let $x_1,\dots,x_n$ denote the Chern roots of the dual
$S^{\vee}$ of the tautological subbundle $S$
of the trivial bundle $\CC^N_{\Gr(n,N)}$ so that
for the total equivariant Chern class of $S$
we have
\ben
c^T(S)=\prod_{i=1}^n(1-x_i).
\een
Then, due to \cite[Lecture~8, Proposition~1.1]{f:ec}
(see also \cite{m:gf}),
the equivariant Schubert classes $\si_{\la}$
can be expressed by
\ben
\si_{\la}=s_{\la}(x|u), \qquad u=(-t_N,\dots,-t_1,0,\dots).
\een
Hence, Theorem~\ref{thm:lrr} yields a multiplication rule
for the equivariant Schubert classes. The corresponding stability
property is implied by Corollary~\ref{cor:finn}.

\bco\label{cor:shub}
We have
\ben
\si_{\la}\ts\si_{\mu}=\sum_{\nu}\ts
d_{\la\mu}^{\ts\nu}\ts\si_{\nu},
\een
where
\beql{lrrulegr}
d_{\la\mu}^{\ts\nu}=\sum_R\sum_{T}
\prod_{\underset{\scriptstyle T(\alpha)
\text{\ts unbarred}}{\alpha\in\la}}
\Big(\ts t^{}_{\tss m+T(\alpha)-c(\alpha)}
-t^{}_{\tss m+T(\alpha)-\rho(\alpha)^{}_{T(\al)}}\ts\Big),
\eeq
summed over all sequences $R$ of the form \eqref{r} and
all $\nu$-bounded
reverse $\la$-tableaux $T\in\Tc(\la,R)$.
In particular, the
$d_{\la\mu}^{\ts\nu}$ are polynomials in the $t_i-t_j$,
$i>j$, with positive integer coefficients.
Moreover, the coefficients $d_{\la\mu}^{\ts\nu}$,
regarded as polynomials in the variables $a_i$
defined in \eqref{at},
are independent of $n$ and $m$, as soon as
the inequalities $n\geqslant \la'_1+\mu'_1$
and $m\geqslant \la_1+\mu_1$ hold.
\qed
\eco

\bex\label{ex:calcschub}
For any $n\geqslant 3$ and $m\geqslant 4$ we have
\ben
\bal
\si_{(2)}\ts \si_{(2,1)}
{}&=\si_{(4,1)}+\si_{(3,2)}+\si_{(3,1,1)}
+\si_{(2,2,1)}\\
{}&+(t_{m+2}-t_{m-1}+t_m-t_{m-2})\ts \si_{(3,1)}
+(t_{m+2}-t_{m-1})\ts \si_{(2,2)}\\
{}&+(t_m-t_{m-1})\ts \si_{(2,1,1)}
+(t_{m+2}-t_{m-1})\ts (t_m-t_{m-1})\ts \si_{(2,1)}.
\eal
\een
This follows from Example~\ref{ex:calc}.
\qed
\eex

The first manifestly positive rule for the
expansion of $\si_{\la}\ts\si_{\mu}$ was given by
Knutson and Tao~\cite{kt:pe} by using combinatorics
of puzzles. Although the stability property
was not pointed out in \cite{kt:pe}, it
can be deduced directly from the puzzle rule or by
applying the weight-preserving bijection
between the puzzles and the
barred tableaux constructed by Kreiman~\cite{k:elr}.

\subsection{Quantum immanants and higher Capelli operators}
\label{subsec:hcap}

Let $\gl_n$ denote the general linear Lie algebra
over $\CC$.
Consider the center $\Z(\gl_n)$ of the universal enveloping algebra
$\U(\gl_n)$. The algebra $\U(\gl_n)$ is equipped with
the natural filtration.
For all $n$ we identify $\gl_{n-1}$ as a subalgebra
of $\gl_n$ in a usual way and denote by $\gl_{\infty}$
the corresponding inductive limit
\ben
\gl_{\infty}=\underset{n}{\bigcup} \ts\ts \gl_n.
\een
Due to Olshanski~\cite{o:ri},
there exist filtration-preserving homomorphisms
\beql{olshh}
o_n:\Z(\gl_n)\to \Z(\gl_{n-1}),\qquad n\geqslant 1,
\eeq
which allow one to define the algebra $\Z$ of the {\it virtual
Casimir elements\/} for the Lie algebra $\gl_{\infty}$
as the inverse limit
\ben
\Z=\lim_{\longleftarrow} \Z(\gl_n),\qquad n\to\infty,
\een
in the category of filtered algebras.

The {\it quantum immanants\/} $\SSb_{\la|n}$ are elements
of the center $\Z(\gl_n)$ of the universal enveloping algebra
$\U(\gl_n)$ parameterized by the diagrams $\la$ with
at most $n$ rows; see \cite{o:qi}.
The elements $\SSb_{\la|n}$ form a basis of $\Z(\gl_n)$
and they are consistent with the Olshanski
homomorphisms \eqref{olshh} so that
\beql{stabqi}
o_n:\SSb_{\la|n}\mapsto \SSb_{\la|n-1},
\eeq
where we assume $\SSb_{\la|n}=0$ if the number of rows of
$\la$ exceeds $n$. For any diagram $\la$, the corresponding
{\it virtual quantum immanant\/}
$\SSb_{\la}$ is then defined as the sequence
\ben
\SSb_{\la}=(\ts\SSb_{\la|n}\ts|\ts n\geqslant 0).
\een
The elements $\SSb_{\la}$ parameterized by all diagrams $\la$
form a basis of the algebra $\Z$
so that we can define the coefficients $f_{\la\mu}^{\ts\nu}$
by the expansion
\ben
\SSb_{\la}\ts \SSb_{\mu}=\sum_{\nu} f_{\la\mu}^{\ts\nu}\ts
\SSb_{\nu}.
\een
Note that the same coefficients
$f_{\la\mu}^{\ts\nu}$ determine the multiplication rule
for the {\it higher Capelli operators\/} $\Delta_{\la}$,
which are defined as the sequences of the images of the
quantum immanants $\SSb_{\la|n}$, where each image is taken
under a natural
representation of $\gl_n$ by differential operators; see
\cite{o:qi, oo:ss}.

\bco\label{cor:posit}
The coefficient $f_{\la\mu}^{\nu}$ is zero
unless $\mu\subseteq\nu$. If $\mu\subseteq\nu$ then
\beql{lrruleqi}
f^{\nu}_{\la\mu}=\sum_R\sum_{T}
\prod_{\underset{\scriptstyle
T(\alpha)\text{\ts unbarred}}{\alpha\in\la}}
\Big(\rho(\alpha)^{}_{T(\al)}-c(\alpha)\Big),
\eeq
summed over all sequences $R$ of the form \eqref{r} and
all $\nu$-bounded
reverse $\la$-tableaux $T\in\Tc(\la,R)$.
In particular, the $f^{\nu}_{\la\mu}$ are nonnegative
integers.
\eco

\bpf
Due to the stability property \eqref{stabqi} of the
quantum immanants, it suffices to calculate the
corresponding coefficients for the expansion of
the products $\SSb_{\la|n}\ts\SSb_{\mu|n}$.
The images of the quantum immanants
$\SSb_{\la|n}$ under the Harish-Chandra
isomorphism can be identified with the double Schur polynomials
$s_{\la}(x\vt a)$ where the sequence $a$ is specialized to
$a_i=-i$; see \cite{o:qi}.
Therefore, the coefficients in question coincide
with the corresponding specializations of
the Littlewood--Richardson polynomials
$c_{\la\mu}^{\nu}(a)$.
\epf

\bex\label{ex:calcqi}
Using Example~\ref{ex:calc} we get
\ben
\SSb_{(2)}\ts \SSb_{(2,1)}
=\SSb_{(4,1)}+\SSb_{(3,2)}+\SSb_{(3,1,1)}
+\SSb_{(2,2,1)}
+5\ts\ts \SSb_{(3,1)}
+3\ts\ts \SSb_{(2,2)}
+\SSb_{(2,1,1)}
+3\ts\ts \SSb_{(2,1)}.
\een
\eex

\medskip

In the course of the proof of
Corollary~\ref{cor:posit}
we also calculated the coefficients
for the expansion of
the products $\SSb_{\la|n}\ts\SSb_{\mu|n}$ for any $n$.
Some other formulas for these coefficients
were obtained in \cite{ms:lr}. In particular, it was shown that
the $f_{\la\mu}^{\ts\nu}$
are integers, although their positivity property
was not established there.

Note also that the
algebra of virtual Casimir elements $\Z$ is isomorphic
to the algebra of {\it shifted symmetric functions\/} $\La^*$;
see \cite{oo:ss}. The latter can be regarded as
the specialization of $\La$ (or rather, its extension over $\CC$)
at $a_i=-i$ for all $i\in\ZZ$.

\section{Supertableau formulas for $s_{\la}(x\vt a)$
and $c_{\la\mu}^{\nu}(a)$}
\label{sec:supt}

Here we obtain one more rule for the calculation
of the Littlewood--Richardson polynomials $c_{\la\mu}^{\nu}(a)$.
It relies on a supertableau representation of the double
Schur polynomials $s_{\la}(x\vt a)$ which is implied
by the results of \cite{m:sf}.
This representation
provides a ``finite" version of the supertableau
formulas of \cite{gg:nt}
and \cite{m:sf}; cf. \cite{cll:fd}.

Fix a positive integer $n$.
For $r\geqslant 1$ set $u^{(r)}=(u_1,\dots,u_r)$
and use the 9th Variation in \cite{m:sf}
with the indeterminates $h_{rs}$ specialized by
\ben
h_{rs}=h_r(u^{(n-r-s+1)})\qquad\text{if}\quad r+s\leqslant n,
\een
and $0$ otherwise, where $h_r$ denotes the $r$-th complete symmetric
polynomial.
Let us write $\wh s_{\la/\mu}(u)$ for the corresponding
Schur functions. Then (8.2) and (9.1) in \cite{m:sf} give
\ben
\wh s_{\la/\mu}(u)=\sum_T \prod_{\alpha\in\la/\mu}
u^{}_{T(\al)},
\een
summed over semistandard tableaux $T$ of shape $\la/\mu$,
such that the entries of the $i$-th row do not exceed $n-\la_i+i$.
Furthermore, using (6.18)\footnote{This formula in \cite{m:sf}
should be corrected by replacing $a^{(\la_j+n-j)}$
with $a^{(\la_i+n-i)}$.}
and (9.6${\tss}'$) in \cite{m:sf} we get
\beql{sdecom}
s_{\la}(x|u)=\sum_{\mu\subseteq\la}
s_{\mu}(x)\ts \wh s_{\la'/\mu'}(-u).
\eeq
Equivalently, this can be interpreted as a combinatorial
expression for the polynomials $s_{\la}(x|u)$ in terms
of ``supertableaux".
Identify the indices of $u$ with the symbols $1',2',\dots$.
A {\it supertableau\/} $T$ is obtained by filling
in the diagram of $\la$ with the indices $1,\dots,n,1',2',\dots$
in such a way that
in each row (resp. column) each primed
index is to the right (resp. below)
of each unprimed index;
unprimed indices weakly increase along the rows and
strictly increase down the columns;
primed indices strictly increase along the rows and
weakly increase down the columns;
primed indices in column $j$ do not exceed $n-\la'_j+j$.
Relation \eqref{sdecom} implies the following.

\bpr\label{prop:supert}
We have
\beql{tab}
s_{\la}(x|u)=
\sum_{T}
\prod_{\underset{\scriptstyle T(\alpha)
\text{\ts\ts unprimed}}{\alpha\in\la}}
x^{}_{T(\al)}
\prod_{\underset{\scriptstyle T(\alpha)
\text{\ts\ts primed}}{\alpha\in\la}} (-u^{}_{T(\al)}),
\eeq
summed over all $\la$-supertableaux $T$.
\qed
\epr

Using \eqref{ssassb}, we get an analogous representation
for the double Schur polynomials $s_{\la}(x\vt a)$.
A {\it reverse supertableau\/} $T$ is obtained by filling
in the diagram of $\la$ with the indices $1,\dots,n,n',(n-1)',\dots$
(including non-positive primed indices)
in such a way that in each row (resp. column) each primed
index is to the right (resp. below)
of each unprimed index;
unprimed indices weakly decrease along the rows and
strictly decrease down the columns;
primed indices strictly decrease along the rows and
weakly decrease down the columns;
primed indices in column $j$ are not less than $\la'_j-j+1$.
The following supertableau representation
of the polynomials $s_{\la}(x\vt a)$
follows from Proposition~\ref{prop:supert}.

\bco\label{prop:superta}
We have
\beql{taba}
s_{\la}(x\vt a)=
\sum_{T}
\prod_{\underset{\scriptstyle T(\alpha)
\text{\ts\ts unprimed}}{\alpha\in\la}}
x^{}_{T(\al)}
\prod_{\underset{\scriptstyle T(\alpha)
\text{\ts\ts primed}}{\alpha\in\la}} (-a^{}_{T(\al)}),
\eeq
summed over all reverse $\la$-supertableaux $T$.
\qed
\eco

\bex\label{ex:supd}
Let $n=2$ and $\la=(2,1)$. By the definition \eqref{defdouble},
\ben
s_{(2,1)}(x\vt a)=(x_2-a_2)(x_1-a_0)(x_1-a_2)+
(x_2-a_2)(x_2-a_1)(x_1-a_2).
\een
On the other hand, the reverse $(2,1)$-supertableaux are

\setlength{\unitlength}{0.75em}
\begin{center}
\begin{picture}(46,4.2)

\put(0,0){\line(0,1){4}}
\put(2,0){\line(0,1){4}}
\put(4,2){\line(0,1){2}}
\put(0,0){\line(1,0){2}}
\put(0,2){\line(1,0){4}}
\put(0,4){\line(1,0){4}}
\put(0.7,0.5){1}
\put(0.7,2.5){2}
\put(2.7,2.5){1}

\put(6,0){\line(0,1){4}}
\put(8,0){\line(0,1){4}}
\put(10,2){\line(0,1){2}}
\put(6,0){\line(1,0){2}}
\put(6,2){\line(1,0){4}}
\put(6,4){\line(1,0){4}}
\put(6.7,0.5){1}
\put(6.7,2.5){2}
\put(8.7,2.5){2}

\put(12,0){\line(0,1){4}}
\put(14,0){\line(0,1){4}}
\put(16,2){\line(0,1){2}}
\put(12,0){\line(1,0){2}}
\put(12,2){\line(1,0){4}}
\put(12,4){\line(1,0){4}}
\put(12.5,0.5){$2^{\tss\prime}$}
\put(12.7,2.5){2}
\put(14.7,2.5){1}

\put(18,0){\line(0,1){4}}
\put(20,0){\line(0,1){4}}
\put(22,2){\line(0,1){2}}
\put(18,0){\line(1,0){2}}
\put(18,2){\line(1,0){4}}
\put(18,4){\line(1,0){4}}
\put(18.5,0.5){$2^{\tss\prime}$}
\put(18.7,2.5){2}
\put(20.7,2.5){2}

\put(24,0){\line(0,1){4}}
\put(26,0){\line(0,1){4}}
\put(28,2){\line(0,1){2}}
\put(24,0){\line(1,0){2}}
\put(24,2){\line(1,0){4}}
\put(24,4){\line(1,0){4}}
\put(24.5,0.5){$2^{\tss\prime}$}
\put(24.7,2.5){1}
\put(26.7,2.5){1}

\put(30,0){\line(0,1){4}}
\put(32,0){\line(0,1){4}}
\put(34,2){\line(0,1){2}}
\put(30,0){\line(1,0){2}}
\put(30,2){\line(1,0){4}}
\put(30,4){\line(1,0){4}}
\put(30.7,0.5){1}
\put(30.7,2.5){2}
\put(32.5,2.5){$0^{\tss\prime}$}

\put(36,0){\line(0,1){4}}
\put(38,0){\line(0,1){4}}
\put(40,2){\line(0,1){2}}
\put(36,0){\line(1,0){2}}
\put(36,2){\line(1,0){4}}
\put(36,4){\line(1,0){4}}
\put(36.7,0.5){1}
\put(36.7,2.5){2}
\put(38.5,2.5){$1^{\tss\prime}$}

\put(42,0){\line(0,1){4}}
\put(44,0){\line(0,1){4}}
\put(46,2){\line(0,1){2}}
\put(42,0){\line(1,0){2}}
\put(42,2){\line(1,0){4}}
\put(42,4){\line(1,0){4}}
\put(42.7,0.5){1}
\put(42.7,2.5){2}
\put(44.5,2.5){$2^{\tss\prime}$}

\end{picture}
\end{center}
\setlength{\unitlength}{1pt}

\setlength{\unitlength}{0.75em}
\begin{center}
\begin{picture}(46,4.2)

\put(0,0){\line(0,1){4}}
\put(2,0){\line(0,1){4}}
\put(4,2){\line(0,1){2}}
\put(0,0){\line(1,0){2}}
\put(0,2){\line(1,0){4}}
\put(0,4){\line(1,0){4}}
\put(0.5,0.5){$2^{\tss\prime}$}
\put(0.7,2.5){2}
\put(2.5,2.5){$0^{\tss\prime}$}

\put(6,0){\line(0,1){4}}
\put(8,0){\line(0,1){4}}
\put(10,2){\line(0,1){2}}
\put(6,0){\line(1,0){2}}
\put(6,2){\line(1,0){4}}
\put(6,4){\line(1,0){4}}
\put(6.5,0.5){$2^{\tss\prime}$}
\put(6.7,2.5){2}
\put(8.5,2.5){$1^{\tss\prime}$}

\put(12,0){\line(0,1){4}}
\put(14,0){\line(0,1){4}}
\put(16,2){\line(0,1){2}}
\put(12,0){\line(1,0){2}}
\put(12,2){\line(1,0){4}}
\put(12,4){\line(1,0){4}}
\put(12.5,0.5){$2^{\tss\prime}$}
\put(12.7,2.5){2}
\put(14.5,2.5){$2^{\tss\prime}$}

\put(18,0){\line(0,1){4}}
\put(20,0){\line(0,1){4}}
\put(22,2){\line(0,1){2}}
\put(18,0){\line(1,0){2}}
\put(18,2){\line(1,0){4}}
\put(18,4){\line(1,0){4}}
\put(18.5,0.5){$2^{\tss\prime}$}
\put(18.7,2.5){1}
\put(20.5,2.5){$0^{\tss\prime}$}

\put(24,0){\line(0,1){4}}
\put(26,0){\line(0,1){4}}
\put(28,2){\line(0,1){2}}
\put(24,0){\line(1,0){2}}
\put(24,2){\line(1,0){4}}
\put(24,4){\line(1,0){4}}
\put(24.5,0.5){$2^{\tss\prime}$}
\put(24.7,2.5){1}
\put(26.5,2.5){$1^{\tss\prime}$}

\put(30,0){\line(0,1){4}}
\put(32,0){\line(0,1){4}}
\put(34,2){\line(0,1){2}}
\put(30,0){\line(1,0){2}}
\put(30,2){\line(1,0){4}}
\put(30,4){\line(1,0){4}}
\put(30.5,0.5){$2^{\tss\prime}$}
\put(30.7,2.5){1}
\put(32.5,2.5){$2^{\tss\prime}$}

\put(36,0){\line(0,1){4}}
\put(38,0){\line(0,1){4}}
\put(40,2){\line(0,1){2}}
\put(36,0){\line(1,0){2}}
\put(36,2){\line(1,0){4}}
\put(36,4){\line(1,0){4}}
\put(36.5,0.5){$2^{\tss\prime}$}
\put(36.5,2.5){$2^{\tss\prime}$}
\put(38.5,2.5){$0^{\tss\prime}$}

\put(42,0){\line(0,1){4}}
\put(44,0){\line(0,1){4}}
\put(46,2){\line(0,1){2}}
\put(42,0){\line(1,0){2}}
\put(42,2){\line(1,0){4}}
\put(42,4){\line(1,0){4}}
\put(42.5,0.5){$2^{\tss\prime}$}
\put(42.5,2.5){$2^{\tss\prime}$}
\put(44.5,2.5){$1^{\tss\prime}$}

\end{picture}
\end{center}
\setlength{\unitlength}{1pt}

\noindent
which yield
\ben
\bal
s_{(2,1)}(x\vt a)&=x_1^2x^{}_2+x_1^{}x^2_2-x^{}_1x^{}_2a^{}_2
-x^2_2a^{}_2-x^2_1a_2-x^{}_1x^{}_2a^{}_0-x^{}_1x^{}_2a^{}_1
-x^{}_1x^{}_2a^{}_2\\
&+x_2^{}a^{}_0a^{}_2+x_2^{}a^{}_1a^{}_2+x_2^{}a^2_2
+x_1^{}a^{}_0a^{}_2+x_1^{}a^{}_1a^{}_2+x_1^{}a^2_2
-a_0^{}a^2_2-a_1^{}a^2_2.
\eal
\een
\eex

Formula \eqref{sdecom} implies
a supertableau representation of
the coefficients $c_{\la\mu}^{\nu}(a,b)$
and hence, of
the Littlewood--Richardson polynomials $c_{\la\mu}^{\nu}(a)$.
The representation for the latter
is neither manifestly positive, nor stable;
it provides an expression for $c_{\la\mu}^{\nu}(a)$
as an alternating sum of monomials in the $a_i$.
Given a sequence $R$ of the form \eqref{r},
construct the set $\Sc(\la,R)$ of barred
reverse $\la$-supertableaux
by analogy with $\Tc(\la,R)$. A tableau $T\in \Sc(\la,R)$
must contain
boxes $\al_1,\dots,\al_l$ occupied by
unprimed indices $r_1,r_2,\dots,r_l$ listed in
the column order which is restricted to the
subtableau of $T$ formed by the unprimed indices.
As before, we distinguish the entries in
$\al_1,\dots,\al_l$ by barring
each of them. For each box $\alpha$ with $\al_i\prec\al\prec\al_{i+1}$,
$0\leqslant i\leqslant l$,
which is occupied by an unprimed index,
set $\rho(\alpha)=\rho^{(i)}$.

\bco\label{cor:supt}
The coefficients $c_{\la\mu}^{\nu}(a,b)$
defined in \eqref{lrcab} can be given by
\beql{lrruleabze}
c^{\nu}_{\la\mu}(a,b)=\sum_R
\sum_{T}
\prod_{\underset{\scriptstyle T(\alpha)
\text{\ts\ts unprimed,\ts unbarred}}{\alpha\in\la}}
a^{}_{\tss T(\alpha)-\rho(\alpha)^{}_{T(\al)}}
\prod_{\underset{\scriptstyle T(\alpha)
\text{\ts\ts primed}}{\alpha\in\la}} (-b^{}_{T(\al)}),
\eeq
summed over sequences $R$ of the form \eqref{r} and
reverse supertableaux $T\in\Sc(\la,R)$.
\eco

\bpf
Applying formula \eqref{sdecom} we can reduce
the calculation of $c^{\nu}_{\la\mu}(a,b)$
to the particular case of the sequence $b=(0)$.
Now \eqref{lrruleabze} follows from Lemma~\ref{lem:msab}.
\epf

\bex\label{ex:superca}
In order to calculate the Littlewood--Richardson polynomial
$c_{(2,1)\ts(2)}^{(2,1)}(a)$ we may take $n=2$;
see Corollary~\ref{cor:finn}. The barred
reverse supertableaux compatible with the sequence $(2)\to (2,1)$
are

\setlength{\unitlength}{0.75em}
\begin{center}
\begin{picture}(34,4.2)

\put(0,0){\line(0,1){4}}
\put(2,0){\line(0,1){4}}
\put(4,2){\line(0,1){2}}
\put(0,0){\line(1,0){2}}
\put(0,2){\line(1,0){4}}
\put(0,4){\line(1,0){4}}
\put(0.7,0.5){1}
\put(0.7,2.5){$\overline{2}$}
\put(2.7,2.5){1}

\put(6,0){\line(0,1){4}}
\put(8,0){\line(0,1){4}}
\put(10,2){\line(0,1){2}}
\put(6,0){\line(1,0){2}}
\put(6,2){\line(1,0){4}}
\put(6,4){\line(1,0){4}}
\put(6.7,0.5){1}
\put(6.7,2.5){$\overline{2}$}
\put(8.7,2.5){2}

\put(12,0){\line(0,1){4}}
\put(14,0){\line(0,1){4}}
\put(16,2){\line(0,1){2}}
\put(12,0){\line(1,0){2}}
\put(12,2){\line(1,0){4}}
\put(12,4){\line(1,0){4}}
\put(12.7,0.5){1}
\put(12.7,2.5){2}
\put(14.7,2.5){$\overline{2}$}

\put(18,0){\line(0,1){4}}
\put(20,0){\line(0,1){4}}
\put(22,2){\line(0,1){2}}
\put(18,0){\line(1,0){2}}
\put(18,2){\line(1,0){4}}
\put(18,4){\line(1,0){4}}
\put(18.5,0.5){$2^{\tss\prime}$}
\put(18.7,2.5){$\overline{2}$}
\put(20.7,2.5){1}

\put(24,0){\line(0,1){4}}
\put(26,0){\line(0,1){4}}
\put(28,2){\line(0,1){2}}
\put(24,0){\line(1,0){2}}
\put(24,2){\line(1,0){4}}
\put(24,4){\line(1,0){4}}
\put(24.5,0.5){$2^{\tss\prime}$}
\put(24.7,2.5){$\overline{2}$}
\put(26.7,2.5){2}

\put(30,0){\line(0,1){4}}
\put(32,0){\line(0,1){4}}
\put(34,2){\line(0,1){2}}
\put(30,0){\line(1,0){2}}
\put(30,2){\line(1,0){4}}
\put(30,4){\line(1,0){4}}
\put(30.5,0.5){$2^{\tss\prime}$}
\put(30.7,2.5){2}
\put(32.5,2.5){$\overline{2}$}

\end{picture}
\end{center}
\setlength{\unitlength}{1pt}

\setlength{\unitlength}{0.75em}
\begin{center}
\begin{picture}(34,4.2)

\put(0,0){\line(0,1){4}}
\put(2,0){\line(0,1){4}}
\put(4,2){\line(0,1){2}}
\put(0,0){\line(1,0){2}}
\put(0,2){\line(1,0){4}}
\put(0,4){\line(1,0){4}}
\put(0.5,0.5){1}
\put(0.7,2.5){$\overline{2}$}
\put(2.5,2.5){$0^{\tss\prime}$}

\put(6,0){\line(0,1){4}}
\put(8,0){\line(0,1){4}}
\put(10,2){\line(0,1){2}}
\put(6,0){\line(1,0){2}}
\put(6,2){\line(1,0){4}}
\put(6,4){\line(1,0){4}}
\put(6.5,0.5){1}
\put(6.7,2.5){$\overline{2}$}
\put(8.5,2.5){$1^{\tss\prime}$}

\put(12,0){\line(0,1){4}}
\put(14,0){\line(0,1){4}}
\put(16,2){\line(0,1){2}}
\put(12,0){\line(1,0){2}}
\put(12,2){\line(1,0){4}}
\put(12,4){\line(1,0){4}}
\put(12.5,0.5){1}
\put(12.7,2.5){$\overline{2}$}
\put(14.5,2.5){$2^{\tss\prime}$}

\put(18,0){\line(0,1){4}}
\put(20,0){\line(0,1){4}}
\put(22,2){\line(0,1){2}}
\put(18,0){\line(1,0){2}}
\put(18,2){\line(1,0){4}}
\put(18,4){\line(1,0){4}}
\put(18.5,0.5){$2^{\tss\prime}$}
\put(18.7,2.5){$\overline{2}$}
\put(20.5,2.5){$0^{\tss\prime}$}

\put(24,0){\line(0,1){4}}
\put(26,0){\line(0,1){4}}
\put(28,2){\line(0,1){2}}
\put(24,0){\line(1,0){2}}
\put(24,2){\line(1,0){4}}
\put(24,4){\line(1,0){4}}
\put(24.5,0.5){$2^{\tss\prime}$}
\put(24.7,2.5){$\overline{2}$}
\put(26.5,2.5){$1^{\tss\prime}$}

\put(30,0){\line(0,1){4}}
\put(32,0){\line(0,1){4}}
\put(34,2){\line(0,1){2}}
\put(30,0){\line(1,0){2}}
\put(30,2){\line(1,0){4}}
\put(30,4){\line(1,0){4}}
\put(30.5,0.5){$2^{\tss\prime}$}
\put(30.7,2.5){$\overline{2}$}
\put(32.5,2.5){$2^{\tss\prime}$}

\end{picture}
\end{center}
\setlength{\unitlength}{1pt}

\noindent
so that
\ben
\bal
c_{(2,1)\ts(2)}^{(2,1)}(a)&=a_{-1}^2+a^{}_{-1}a^{}_1+a^{}_{-1}a^{}_2
-a^{}_{-1}a^{}_2-a^{}_{1}a^{}_2-a^2_2
\\
{}&-a^{}_{-1}a^{}_0-a^{}_{-1}a^{}_1-a^{}_{-1}a^{}_2
+a^{}_{0}a^{}_2+a^{}_{1}a^{}_2+a^2_2\\
{}&=a_{-1}^2-a^{}_{-1}a^{}_0
-a^{}_{-1}a^{}_2+a^{}_{0}a^{}_2,
\eal
\een
which agrees with Example~\ref{ex:calc}.
\qed
\eex

\end{document}